\documentclass{amsart}[12pt]
\usepackage{amsmath,amssymb,enumerate}

\usepackage{color}

\theoremstyle{plain}
\newtheorem{thm}{Theorem}[section]
\newtheorem{lem}[thm]{Lemma}
\theoremstyle{definition}
\newtheorem{defi}[thm]{Definition}

\newtheorem{remark}[thm]{Remark}

\newtheorem{prop}[thm]{Proposition}
\newtoks\by
\newtoks\paper
\newtoks\book
\newtoks\jour
\newtoks\yr
\newtoks\pages
\newtoks\vol
\newtoks\publ
\def\ota{{\hbox\vol{???}}}
\def\cLear{\by=\ota\paper=\ota\book=\ota\jour=\ota\yr=\ota
\pages=\ota\vol=\ota\publ=\ota}
\def\endpaper{\the\by, \the\paper.
{\it\the\jour\/} {\bf \the\vol} (\the\yr), \the\pages.\cLear}
\def\endbook{\the\by, {\it\the\book}. \the\publ.\cLear}
\def\endprep{\the\by, \the\paper. \the\jour.\cLear}
\def\name#1#2{#1 #2}
\def\et{ and }

\def\supp{\textup{supp }}

\def\id{\textup{id}}
\def\spa{\textup{span}}
\def\sign{\textup{sign}}
\def\be{\begin{equation}}
\def\ee{\end{equation}}

\numberwithin{equation}{section} \headheight=12pt
\begin{document}
\title[Embeddings between Lorenz \dots]{Embeddings between Lorenz sequence spaces are strictly singular}

\author[J.Lang]{J. Lang }
\address{Department of Mathematics, The Ohio State University,
231 West 18th Avenue, Columbus, OH 43210, USA}
\email{lang@math.osu.edu}
\author[A.Nekvinda]{A. Nekvinda}
\address{Department of Mathematics\\ Faculty of Civil Engineereng\\
       Czech Technical University\\Th\' akurova 7\\16629 Prague 6\\
       Czech Republic}
\email{nales@mat.fsv.cvut.cz}
\date{}
\thanks{The second author was supported by the European Regional Development Fund,
project No. CZ $02.1.01/0.0/0.0/16-019/0000778$}

\maketitle
\begin{abstract}
Given $0<p,q, r<\infty $ and $ q<r\le \infty$ we consider the natural embedding $\ell_{p,q}\hookrightarrow \ell_{p,r}$ between Lorenz sequence spaces. We prove that this non-compact embedding is always strictly singular but not finitely strictly singular.
\end{abstract}
\section{Introduction}

Let $T:X\to Y$ be a linear operator between Banach spaces, we recall that $T$ is a strictly singular  if there is no infinite dimensional closed subspace $Z$ of $X$ such that $T:Z \to T(Z)$, the restriction of $T$ to $Z$, is an isomorphism.	More over we say that $T$ is finitely strictly singular if: for every $\varepsilon >0$, there exists $n_{\varepsilon} \ge 1$ such that every subspace $E$ of $X$ with dimension greater that $n_{\varepsilon}$, there exists $x$ in the unit sphere of $E$ such that $\|T(x)\|_Y \le \varepsilon.$

Strictly singular operators and finitely strictly singular operators, which encompass compact operators, possess same  properties which are habitually connected with compact operators. For example it is well know that Fredholm operators are invariant when  perturbed by strictly singular operators
(i.e. if $T$ is Fredholm and $S$ is strictly singular then $T+S$ is Fredholm, see \cite[Therem 4.63]{AA book}). 

Many of non-compact operators in Analysis are strictly singular or finitely strictly singular.
For example Fourier transform, which is obviously non-compact, when is considered as a map from $L^p$ into $L^{p'}$, is finitely strictly singular for $1< p<2$ and strictly singular when $p=1$  (see \cite{LRP}). 

The natural embedding of sequence spaces
\[ I:l^p \to l^q, \quad \mbox{ for } p<q,
\] which is non-compact is finitely strictly singular (see \cite{Pl}).

Also the limiting optimal Sobolev embedding $E_d$ into continuous functions
\begin{align} \label{eq E_d}
	&E_d:W^{1}_{0}L^{d,1} ((0,1)^d )\hookrightarrow C((0,1)^d ),
\end{align}
(where $W^{1}_{0}L^{d,1} ((0,1)^d )$ denote a space of all functions $u$ for which $|\nabla u| $ belongs to Lorentz space $L^{d,1}$  and $u$ has a zero trace), is non-compact but finitely strictly singular (see \cite{LM}).

From the above examples arises a natural and quite intrigue question: Whether all limiting Sobolev embeddings on bounded smooth domain are strictly singular or finitely strictly singular?

Since the optimal target spaces for Sobolev embeddings are quite often  Lorentz spaces then in order to investigate the above question, one needs to get more information about natural embedding between Lorentz sequence spaces:
\[I: l_{p,q} \to l_{r,s}. 
\]  

This paper continues an investigation started in \cite{LN} where was proved that the natural embedding $l_{p} \to l_{p,\infty}$ for $p>1$ is strictly singular.

 The paper is structured as follows. In Sect. 2, we recall the definitions we use,
 and we collect all necessary later-needed material and technical lemmas. In Sect. 3 is proved that embedding $l_{p,q} \to l_{p,r}$ is strictly singular for $q<r\le\infty$ and in Sect. 4  we showed that this embedding  is not finitely strictly singular.

\section{preliminaries}

In this section we recall definitions, notations and some technical lemmas needed in Sections 3 and 4.
We start by recalling definition  of strictly singular and finitely strictly singular operators.
\begin{defi}
A bounded operator $T:X \to Y$ between Banach spaces is said to be strictly singular if there is no infinite dimensional closed subspace $Z$ of $X$ such that $T:Z \to T(Z)$, the restriction of $T$ to $Z$, is an isomorphism.	
\end{defi}
	
See \cite[section 4.5]{AA book} for more about strictly singular operators.

\begin{defi}
	An operator $T$ from a Banach space $X$ into a Banach space $Y$ is finitely strictly singular if: for every $\varepsilon >0$, there exists $n_{\varepsilon} \ge 1$ such that every subspace $E$ of $X$ with dimension greater that $n_{\varepsilon}$, there exists $x$ in the unit sphere of $E$ such that $\|T(x)\|_Y \le \varepsilon.$
\end{defi}

Let $T:X\to Y$ be a linear map between Banach spaces then the $n$-th  Bernstein numbers (or Bernstein widths) is defined by
	$$
	b_n(T)=\sup_{ E\subset X,\dim(E)=n } \inf_{f\in E,\|f\|_X=1}  \|T(f)\|_{Y}.
	$$
It is not too hard to see that the operator $T$ is finitely strictly singular if and only if $b_n(T) \to 0$ and that we have the following relations:
\[ \mbox{compact} \quad \Rightarrow  \quad \mbox{finitely strictly singular}  \quad \Rightarrow  \quad \mbox{strictly singular.}\]

We consider in this paper a little more general concept of quasi-Banach spaces which satisfy the "triangle" inequality with a constant.
Denote for $u=(u(1),u(2),\dots)$ the modulus sequence $|u|=(|u(1)|, |u(2)|, \dots)$. We say that $|u|\le |v|$ if $|u(i)|\le |v(i)|$ for each $i\in\mathbb{N}$.
\begin{defi}
Let $\mathcal{S}$ be a a set of all sequences of real numbers and $\|.\|:\mathcal{S}\rightarrow[0,\infty]$. Assume that $\|.\|$ satisfies  for all $u,v\in \mathcal{S}$ and $\alpha\in\mathbb{R}$ we have
\begin{enumerate}[\rm(i)]
\item  $\|u+v\|\le T(\|u\|+\|v\|)$ for some $T\ge 1$,
\item $\|\alpha u\|=|\alpha|\ \|u\|$,
\item $\|u\|\ge 0$ and $\|u\|= 0$ if and only if $u=0$,
\item $\|u\|=\|\ |u|\ \|$,
\item if $|u|\le |v|$ then  $\|u\|\le \|v\|$,
\item if $0\le u_n\nearrow u$ then $\|u_n\|\nearrow\|u\|$,
\item if $\#\{i;u(i)\neq 0\}<\infty$ then $\|u\|<\infty$.
\end{enumerate}
Define $X:=\{u;\|u\|<\infty\}$. Then we call $X$ a sequence quasi-Banach function space.
\end{defi}
By an analogous way we could define a quasi-Banach function space of functions on a domain $\Omega$.
Remark that each quasi-Banach function space is complete (for details see for instance \cite{NP}, Corollary 3.7).

We can extend the definition of strictly singular operators on quasi-Banach spaces by the following alternative definition:
\begin{defi}
Let $X,Y$ be quasi-Banach spaces and assume that $T:X\rightarrow Y$ be a linear bounded operator. We say that $T$ is strictly singular operator if
\begin{align*}
&\inf\{\|Tx\|_Y;\|x\|_X=1, x\in Z\}=0
\end{align*}
for each infinite dimensional subspace $Z\subset X$.
\end{defi}

For a finite set $F$ denote by $\#(F)$ the number of elements of $F$.

\begin{defi}
Given a sequence $a=(a(1),a(2),\dots)\in c_0$ we set for $\lambda>0$
\begin{align*}
&\mu_a(\lambda)=\#\{i;|a(i)|>\lambda\}
\end{align*}
and
\begin{align*}
&a^*(j)=\min\{\lambda>0;\mu_a(\lambda)\le j\}.
\end{align*}
Define $a^*=(a^*(1),a^*(2),\dots)$ a non-increasing rearrangement of $a$.
\end{defi}
For a sequence $a=(a(1),a(2),\dots)\in c_0$ denote
\begin{align*}
&\supp a=\{j\in \mathbb{N}; a(j)\neq 0\}.
\end{align*}

\begin{defi}
Given a sequence $u=(u(1),u(2),\dots)\in c_0$ with $\supp u=\{n_1,n_2,...,n_k\}\subset\mathbb{N}$ and $n_1<n_2<\dots n_k$. Define a non-increasing rearrangement $u^{\diamond}$ of $u$ with respect to $\supp u$ by
\begin{align*}
&\begin{cases}
u^{\diamond}(n_j)=u^*(j)&\ j\in\{1,2,\dots,k\},\\
u^{\diamond}(i)=0&\ i\notin\{n_1,n_2,...,n_k\}.
\end{cases}
\end{align*}

\end{defi}
\begin{remark}
If $\supp u:=\{n+1,n+2,\dots,m\}$ then
\begin{align}
&u^{\diamond}(j)=u^*(j-n).\label{vjrvjvjvvvj}
\end{align}
\end{remark}

In the next we recall the definition of sequence Lorentz spaces.
\begin{defi}
Let $p\in (0,\infty),q\in (0,\infty]$. Define for a sequence $u$
\begin{align*}
&\|u\|_{p,q}=
\begin{cases}
\Big(\sum\limits_{j=1}^\infty j^{q/p-1}(u^*(j))^q \Big)^{1/q}&\ \ \text{if}\ \ q<\infty,\\
\sup\{j^{1/p}\ u^*(j); j=1,2\dots\}&\ \ \text{if}\ \ q=\infty.
\end{cases}
\end{align*}
We define Lorentz space $l_{p,q}$ as a collection of all sequences $u$ for which the norm $\|u\|_{p,q}$ is finite.
\end{defi}

Given $u\in \ell_{p,q}$ we will write $u(i)$ for the value of $u$ at the index $i$.

\begin{lem}
Let $0<p<\infty, 0<q\le \infty$.
The space $l_{p,q}$ is a quasi-Banach function space.
\end{lem}
\begin{proof}
As in \cite{BS} (see (1.16) in Proposition 1.7) we can prove
\begin{align*}
&(u+v)^*(i+j)\le u^*(i)+v^*(j).
\end{align*}
Split the sum
\begin{align*}
&\sum\limits_{j=1}^\infty j^{q/p-1}(u+v)^*(j))^q
\end{align*}
into two sums, the first one is  over odd numbers, the second one is over even numbers. For both sums we can easily prove the quasi-triangle inequality. The other properties are easy.
\end{proof}

\begin{lem}\label{vjifdjkvfjo}
Let $0< p<\infty,0< q<\infty$. Then we have for all $n\in\mathbb{N}$
\begin{align*}
&\Big(\sum_{j=1}^n j^{q/p-1}\Big)^{1/q}\approx n^{1/p}.
\end{align*}
\end{lem}
\begin{proof}
For all $n$ we have
\begin{align*}
&\sum_{j=1}^n j^{q/p-1}\approx \int_0^{n} t^{q/p-1}\ dt=\frac{p}{q}n^{q/p}\approx  n^{q/p}.
\end{align*}
\end{proof}

\begin{prop}
Let $0< p<\infty, 0< q<r\le\infty$. Then  $\ell^{p,q}\hookrightarrow \ell^{p,r}$. Denote by $D_{p,q}$ the norm of this embedding, i.e.
\begin{align}
&\|a\|_{p,r}\le D_{q,r}\|a\|_{p,q}\label{dcehoiefroir}
\end{align}
for all sequences $a$.
\end{prop}

\begin{defi}
	Let $X$ be a quasi-Banach function space of functions defined over $\Omega$. We say that $f\in X$ has an absolutely continuous norm in $X$, written $f\in X_a$, if for every non-increasing sequence of measurable sets $G_n \subset \Omega$ with $|G_n | \searrow 0$ we have $\|f \chi_{G_n} \| \searrow 0$. We say that $X$ has an absolutely continuous norm if $X_a=X$.
	
\end{defi}

\begin{lem}
Let $0< p<\infty,0< q<\infty$. Then $\ell_{p,q}$ has an absolutely continuous norm.
\end{lem}
\begin{proof}
Take $u\in\ell_{p,q}$. Set
\begin{align*}
&u_n(j)=\begin{cases}
0&1\le j\le n,\\
u(j)&n+1\le j.
\end{cases}
\end{align*}
Since $\|u\|_{p,q}\le K<\infty$ we have by \eqref{dcehoiefroir} for each $n$
\begin{align*}
&K\ge\Big(\sum\limits_{j=1}^\infty j^{q/p-1}(u^*(j))^q \Big)^{1/q}\gtrsim n^{1/p}u^*(n)
\end{align*}
and so
\begin{align*}
&u^*(n)\lesssim n^{-1/p}.
\end{align*}
It implies for any $j\in\mathbb{N}$ that $\lim_{n\rightarrow\infty}u_n^*(j)= 0$ and consequently, due to the Lebesgue dominated convergence theorem we obtain
\begin{align*}
&\|u_n\|_{p,q}=\Big(\sum\limits_{j=1}^\infty j^{q/p-1}(u_n^*(j))^q \Big)^{1/q}\rightarrow 0\ \ \text{for}\ \ n\rightarrow\infty.
\end{align*}
\end{proof}

For a sequence $b=(b(1),b(2),\dots)$ and $m\in\mathbb{N}$
set
\begin{align*}
	&P_m(b)=(b(1),b(2),\dots,b(m),0,0,\dots)\\
	&R_m(b)=b-P_m b=(0,0,\dots,0,b(m+1),b(m+2),\dots).
\end{align*}

Let $X\subset \ell_{p,q}$ be a closed subspace with $\dim X=\infty$. Define $X_m=R_m (X)$. It is easy to see that $X_m$ is a closed subspace with $\dim X_m=\infty$.

	Let $0< p<\infty,0< q\le \infty$. Since $\ell_{p,q}$ is a sequence Banach function space we have $T\ge 1$ such that
	\begin{align*}
		&\|u+v\|_{p,q}\le T(\|u\|_{p,q}+\|u\|_{p,q}).
	\end{align*}

Remark that it implies directly
\begin{align}
	&u=v+w\ \Rightarrow\ \|v\|_{p,q}\ge \frac{1}{T}\|u\|_{p,q}-\|w\|_{p,q},\label{sdcdiovnhdo}\\
	&\Big\|\sum_{j=1}^n u_j\Big\|_{p,q}\le \sum_{j=1}^n T^j\|u_j\|_{p,q}.\label{dovdovio}
\end{align}

\begin{lem}\label{sdcdlvhndvnhvki}
Let $0< p<\infty,0< q<\infty$ and $\alpha>0$. Assume $v_j\in\ell_{p,q}$ have pairwise disjoint supports and $\|v_j\|_{p,q}\ge \alpha$.
Then
\begin{align*}
&\lim_{k\rightarrow\infty}\Big\|\sum_{j=1}^k v_j\Big\|_{p,q}= \infty.
\end{align*}
\end{lem}
\begin{proof}
Since $v_j$ have pairwise disjoint supports we can write
\begin{align*}
&\Big\|\sum_{j=1}^k v_j\Big\|_{p,q}=\Big\|\sum_{j=1}^k |v_j|\Big\|_{p,q}.
\end{align*}
Assume that there exists a positive constant $C$ independent of $k$ such that
\begin{align*}
&C\ge \Big\|\sum_{j=1}^k v_j\Big\|_{p,q}=\Big\|\sum_{j=1}^k |v_j|\Big\|_{p,q}.
\end{align*}
Since
\begin{align*}
&\Big\|\sum_{j=1}^k |v_j|\Big\|_{p,q}\nearrow \Big\|\sum_{j=1}^\infty |v_j|\Big\|_{p,q}
\end{align*}
we have
\begin{align*}
&C\ge \Big\|\sum_{j=1}^\infty |v_j|\Big\|_{p,q}.
\end{align*}
By the absolute continuity of $\|.\|_{p,q}$ we obtain
\begin{align*}
&\alpha\le \|v_n\|_{p,q}=\|\ |v_n|\ \|_{p,q}\le \Big\|\sum_{j=n}^\infty |v_j|\Big\|_{p,q}\rightarrow 0
\end{align*}
which is a contradiction.
\end{proof}

\begin{lem}\label{adkcvhovfj}
Suppose $0< p<\infty,0< q<\infty$. Let $X\subset \ell_{p,q}$ be a closed subspace with $\dim X=\infty$.
Assume $n, N\in\mathbb{N}$ and $\varepsilon>0,\frac{1}{T}> \delta>0$. Then there exists $m\in\mathbb{N}$ and $u\in X_n$ such that denoting $v:=P_m u$, $w:=R_m u$
\begin{align}
&\|u\|_{p,q}=1,\label{evklfkobb}\\
&m> 2n,\ m\ge N\label{dcvwvjwvjj}\\
&\supp v\subset\{n+1,n+2,\dots,m\},\label{dlkdlkvdkvk}\\
&|v(j)|\le \varepsilon\text{ for all }j,\label{skcnlvwig}\\
&\frac{1}{T}-\delta\le \|v\|_{p,q}\le1,\label{odcjvjevjvj}\\
&\|w\|_{p,q}\le \delta.\label{spvdjpbjpojb}
\end{align}
\end{lem}
\begin{proof}Set $n_0:=n$ and
construct by induction sequences $n_0< n_1<n_2<\dots$ and $u_i\in X$ such that setting $v_i:=P_{n_{i}}u_i$, $w_i:=R_{n_{i}}u_i$ we have
\begin{align}
&\supp v_i\subset\{n_{i-1}+1,n_{i-1}+2,\dots,n_{i}\},\label{dkdcdjpvdpvv}\\
&\frac{1}{T}-\frac{\delta}{(2T)^{i}}\le \|v_i\|_{p,q}\le 1.\label{wdkvjwpvjwvppj}\\
&\|w_i\|_{p,q}\le \frac{\delta}{(2T)^{i}}.\label{dlkvdnvkdjvdkv}
\end{align}

Since $\dim X_n=\infty$ we can find $u_1\in X_n$ with $\|u_1\|_{p,q}=1$. Take $n_1>n$ such that $\|R_{n_1}u_1\|_{p,q}\le \delta/(2T)$. Denote $v_1:=P_{n_1}u_1$, $w_1:=R_{n_1}u_1$. Clearly, $\supp v_1\subset\{n+1,n+2,\dots,n_1\}$ and
\begin{align*}
&1\ge \|v_1\|_{p,q}\overset{\eqref{sdcdiovnhdo}}{\ge} \frac{1}{T}\|u_1\|_{p,q}-\|w_1\|_{p,q}\ge \frac{1}{T}-\frac{\delta}{2T}.
\end{align*}

Suppose that we have constructed $n_0< n_1<n_2<\dots<n_k$, $u_1,u_2,\dots,u_k\in X$ and appropriate functions $v_1,v_2,\dots,v_k$ satisfying  \eqref{dkdcdjpvdpvv} and \eqref{wdkvjwpvjwvppj}.
Since $\dim X_{n_k}=\infty$ we are able to find $u_{k+1}\in X_{n_k}$ with $\|u_{k+1}\|_{p,q}=1$. It is easy to see that we can take an index $n_{k+1}\ge n_k$ such that
$\|R_{n_{k+1}}u_{k+1}\|_{p,q}\le \frac{\delta}{(2T)^{k+1}}$. Set $w_{k+1}=R_{n_{k+1}}u_{k+1}$, $v_{k+1}=P_{n_{k+1}}u_{k+1}$. Consequently
\begin{align*}
&1\ge \|v_{k+1}\|_{p,q}\overset{\eqref{sdcdiovnhdo}}{\ge} \frac{1}{T}\|u_{k+1}\|_{p,q}-\|w_{k+1}\|_{p,q}\ge \frac{1}{T}-\frac{\delta}{(2T)^{k+1}}.
\end{align*}
Moreover $\supp v_{k+1}\subset\{n_{k}+1,n_{k}+2,\dots,n_{k+1}\}$.

Consider now sequences
\begin{align*}
&y_k:=\sum_{j=1}^k u_j,\ \  \ s_k:=\|y_k\|_{p,q}.
\end{align*}
By \eqref{sdcdiovnhdo}, \eqref{dovdovio} and \eqref{dlkvdnvkdjvdkv} we can write
\begin{align*}
&s_k=\Big\|\sum_{j=1}^k u_j\Big\|_{p,q}=\Big\|\sum_{j=1}^k v_j+\sum_{j=1}^k w_j\Big\|_{p,q}\overset{\eqref{sdcdiovnhdo}}{\ge} \frac{1}{T}\Big\|\sum_{j=1}^k v_j\Big\|_{p,q}-\Big\|\sum_{j=1}^k w_j\Big\|_{p,q}\\
&\overset{\eqref{dovdovio}}{\ge} \frac{1}{T}\Big\|\sum_{j=1}^k v_j\Big\|_{p,q}-\sum_{j=1}^k T^j\|w_j\|_{p,q}\overset{\eqref{dlkvdnvkdjvdkv}}{\ge} \frac{1}{T}\Big\|\sum_{j=1}^k v_j\Big\|_{p,q}-\sum_{j=1}^k T^j\frac{\delta}{(2T)^{j}}\\
&\ge \frac{1}{T}\Big\|\sum_{j=1}^k v_j\Big\|_{p,q}-\delta.
\end{align*}
Since by \eqref{wdkvjwpvjwvppj} we obtain
\begin{align*}
& \|v_i\|_{p,q}\ge\frac{1}{T}-\frac{\delta}{(2T)^{i}}\ge\frac{1}{T}-\frac{\delta}{(2T)}:=\alpha>0
\end{align*}
and $v_j$ have pairwise disjoint supports by \eqref{dlkdlkvdkvk}, Lemma \ref{sdcdlvhndvnhvki} gives
\begin{align*}
&\lim_{k\rightarrow\infty}\Big\|\sum_{j=1}^k v_j\Big\|_{p,q}=\infty
\end{align*}
and consequently $s_k\nearrow\infty$.

Then we are able to find $m$ large enough such that
\begin{align}
&m\ge 2n,\ m\ge N,\ \frac{1}{s_m}\le \varepsilon\label{wdvfhnrogronhi}
\end{align}
and set
\begin{align*}
&u=\frac{1}{s_m}\sum_{j=1}^m u_j.
\end{align*}

It is seen from the definition of $s_m$
\begin{align*}
&\|u\|_{p,q}= 1
\end{align*}
which proves \eqref{evklfkobb}.

Clearly, condition \eqref{dcvwvjwvjj} is satisfied. By the definition of $v=P_m u$ we obtain directly
\begin{align*}
&\supp v\subset\{n_0+1,n_0+2,\dots,m\}=\{n+1,n+2,\dots,m\}
\end{align*}
which proves \eqref{dlkdlkvdkvk}.

Fix now $j\in\mathbb{N}$. If $j> m$ we have $|v(j)|=0$.

Assume $j\le m$. Assume $\|u_k\|_{p,q}=1$. If there is $j\in\mathbb{N}$ with $|u_k(j)|>1$ then we have immediately $\|u_k\|_{p,q}>1$. Thus we have for each $s\in\mathbb{N}$
\begin{align*}
&|u_k(s)|\le 1.
\end{align*}
Clearly, using that $u_j$ have pairwise disjoint supports, we have for each $s$
\begin{align*}
&|v(s)|\le |u(s)|\le \frac{1}{s_m}\sum_{j=1}^m |u_j(s)|\lesssim \frac{1}{s_m}\overset{\eqref{wdvfhnrogronhi}}{\le} \varepsilon
\end{align*}
which proves \eqref{skcnlvwig}.

At last,
\begin{align*}
&w=R_mu=R_m\Big(\frac{1}{s_m}\sum_{j=1}^m u_j\Big)=R_m\Big(\frac{1}{s_m}\sum_{j=1}^m v_j+\frac{1}{s_m}\sum_{j=1}^m w_j\Big)\\
&\frac{1}{s_m}\sum_{j=1}^m R_m( v_j)+\frac{1}{s_m}\sum_{j=1}^m R_m( w_j)=\frac{1}{s_m}\sum_{j=1}^m R_m( w_j).
\end{align*}
Thus
\begin{align*}
&\|w\|_{p,q}\overset{\eqref{dovdovio}}{\le}\frac{1}{s_m}\sum_{j=1}^m T^j \|R_m( w_j)\|_{p,q}\le\frac{1}{s_m}\sum_{j=1}^m T^j \|w_j\|_{p,q}\\
&\overset{\eqref{dlkvdnvkdjvdkv}}{\le}\frac{1}{s_m}\sum_{j=1}^m T^j \frac{\delta}{(2T)^{i}}\le\frac{1}{s_m}\sum_{j=1}^m \frac{\delta}{2^j}\le\frac{\delta}{s_m}\le \delta
\end{align*}
which proves \eqref{spvdjpbjpojb}.

Finally,  The property \eqref{odcjvjevjvj} follows directly from
\begin{align*}
&1\ge \|v\|_{p,q}\overset{\eqref{sdcdiovnhdo}}{\ge}\frac{1}{T}\|u\|_{p,q}-\|w\|_{p,q}\ge\frac{1}{T}-\delta
\end{align*}
which finishes the proof.
\end{proof}

\section{The embedding is strictly singular}
 We prove in this section that the embedding $\ell_{p,q}\hookrightarrow\ell_{p,r}$ is strictly singular for $0< p<\infty$, $q<r\le \infty$. First we prove this assertion for $q<\infty$.

\begin{thm}\label{wdjkchdocvho}
Let $0< p, q,r<\infty$, $q<r$. Then the embedding $\ell_{p,q}\hookrightarrow\ell_{p,r}$ is strictly singular.
\end{thm}
\begin{proof}

Let $X\subset \ell_{p,q}$ be a closed subspace with $\dim X=\infty$ and fix a sequence $\widetilde{a}\in \ell^{p,q}$, $\widetilde{a}(1)\ge\widetilde{a}(2)\ge\dots>0$ such that
\begin{align}
&0<\|\widetilde{a}\|_{p,q}\le 1.\label{flkkbjbb}
\end{align}

Having a sequence $0=n_0<n_1<n_2<\dots$ and $u_k\in X_{n_{k-1}}=R_{n_{k-1}}(X)$ we denote for $k\ge 1$
\begin{align*}
&v_k=P_{n_k} u_k,\ w_k=R_{n_k} u_k,\\
&I_k=\{n_{k-1}+1,n_{k-1}+2,\dots,n_k\},\\
&b_k=\min\{|v_k(j)|;v_{k}(j)\neq 0,j\in I_k\}.
\end{align*}

Choose $0<\delta<1/T$.

We will construct by mathematical induction a sequence of integers $0=n_0<n_1<n_2<\dots$, a sequence of positive real numbers $\varepsilon_1,\varepsilon_2,\dots$, a sequence of functions $u_k\in X_{n_{k-1}}$, $k\ge 1$, and a fixed sequence $a(1)\ge a(2)\ge\dots>0$ with the following properties. We set
\begin{align}
&c_k=\min\Big\{\frac{1}{k n_{k}j^{q/p-1}};j=1,2\dots,n_{k}\Big\}\label{fvdjvjfvfnh}
\end{align}
and we have for $k\ge 1$
\begin{align}
&\|u_k\|_{p,q}= 1,\label{bkpjpnjn}\\
&2n_{k-1}<  n_k.\label{elfkbekbenk}\\
&\supp v_k\subset I_k,\label{bjklnkbkpb}\\
&\varepsilon_{k+1}\le \min\{b_{k},c_{k}^{1/q},a(n_{k})\}\label{scjkvhsdjovh}\\
&\varepsilon_{k+1} n_{k}^{1/p}\le 1,\label{efbklbnpbkk}\\
&|v_{k}(j)|\le \varepsilon_k\ \text{for}\ j\in\mathbb{N},\label{ebjklnbkbnkpn}\\
&a(j)\le \widetilde{a}(j)\ \text{for}\ j\in\mathbb{N},\label{skjcdwlkjcjdwc}\\
&a(n_k+1)\le b_k,\label{dsvjdvjvjwvkjv}\\
&\frac{1}{T}-\delta\le \|v_k\|_{p,q}\le1,\label{efbklbpmjkpbn}\\
&\|w_k\|_{p,q}\le \frac{\delta}{(2T)^k}.\label{elbkjbpbjp}
\end{align}

Consider first $k=1$. Find $u_1\in X_0=X$ with $\|u_1\|_{p,q}=1$ a nd set $\varepsilon_1=1$. There exists $n_1>n_0$ such that $\|w_1\|_{p,q}\le \delta/(2T)$ and set $a(i)=\widetilde{a}(i)$, $i\in I_1$. Clearly,
\begin{align}
&1=\|u_1\|_{p,q}\ge\|v_1\|_{p,q}\overset{\eqref{sdcdiovnhdo}}{\ge}\frac{1}{T}\|u_1\|_{p,q}-\|w_1\|_{p,q}\ge \frac{1}{T}-\delta.
\end{align}
Now, it is easy to verify conditions \eqref{bkpjpnjn} -- \eqref{elbkjbpbjp}.

Suppose that we have constructed $n_0< n_1<n_2<\dots<n_k$, $\varepsilon_i,c_i$ for $1\le i\le k$, the sequence $a(i)$ for $i\in I_1\cup I_2\cup\dots\cup I_k$ and functions $u_1,u_2,\dots,u_k\in X$ satisfying the above conditions.

Choose $\varepsilon_{k+1}$ such that
\begin{align}
&\varepsilon_{k+1}\le \min\{b_{k},c_{k}^{1/q},a(n_{k})\}, \ \ \ \varepsilon_{k+1} n_{k}^{1/p}\le 1. \label{djvchdwovhwvoih}
\end{align}
According to Lemma \ref{adkcvhovfj} with $\varepsilon:=\varepsilon_{k+1}$ and $\delta:=\frac{\delta}{(2T)^{k+1}}$ there exists $m>2n_k$ and $u$ with $\|u\|_{p,q}=1$ such that \eqref{evklfkobb}-\eqref{spvdjpbjpojb} are satisfied with $v:=P_m u$, $w:=R_m u$. Set
\begin{align*}
&n_{k+1}:=m,\ \ \ u_{k+1}:=u.
\end{align*}
Then
\begin{align*}
&v=v_{k+1}=P_{n_{k+}}u_{k+1},\ \ \ w=w_{k+1}=R_{n_{k+}}u_{k+1}.
\end{align*}
Set
\begin{align}
&\lambda_{k}:=\min\Big\{\frac{a(n_k)}{\widetilde{a}(n_k+1)},\frac{b_k}{\widetilde{a}(n_k+1)},1\Big\} \label{wdweijefjf}
\end{align}
and
\begin{align}
&a(j):=\lambda_{k}\widetilde{a}(j),\ \ j\in I_{k+1}.\label{kdljvcpwvjpewvj}
\end{align}

Now, \eqref{elbkjbpbjp} follows from \eqref{spvdjpbjpojb}.

Further
\begin{align*}
&1\ge \|v_{k+1}\|_{p,q}\overset{\eqref{odcjvjevjvj}}{\ge} \frac{1}{T}-\frac{\delta}{(2T)^{k+1}}\ge \frac{1}{T}-\delta
\end{align*}
which proves \eqref{efbklbpmjkpbn}.

The properties \eqref{scjkvhsdjovh} and \eqref{efbklbnpbkk} are an immediate consequence of choosing of $\varepsilon_{k+1}$ which is done in \eqref{djvchdwovhwvoih}.

The property \eqref{ebjklnbkbnkpn} follows directly from \eqref{skcnlvwig}. Moreover, by \eqref{wdweijefjf} and \eqref{kdljvcpwvjpewvj} we obtain
\begin{align*}
&a(n_k+1)=\lambda_k\widetilde{a}(n_k+1)\le b_k
\end{align*}
which confirms \eqref{dsvjdvjvjwvkjv}.

Verify that $a(i)$ is non-increasing. If $i,j\in I_{k+1}$ then
\begin{align*}
&a(i)=\lambda_k\widetilde{a}(i)\ge \lambda_k\widetilde{a}(j)=a(j).
\end{align*}
Moreover
\begin{align*}
&a(n_k+1)=\lambda_k\widetilde{a}(n_k+1)\overset{\eqref{wdweijefjf}}{\le}\frac{a(n_k)}{\widetilde{a}(n_k+1)}\widetilde{a}(n_k+1)=a(n_k)
\end{align*}
and $a(i)$ is really non-increasing.

By \eqref{wdweijefjf} we have $\lambda_k\le 1$ and so by \eqref{kdljvcpwvjpewvj} we have \eqref{skjcdwlkjcjdwc}.

At last, properties \eqref{dlkdlkvdkvk}, \eqref{dcvwvjwvjj} and \eqref{evklfkobb} give properties \eqref{bjklnkbkpb}, \eqref{elfkbekbenk} and \eqref{bkpjpnjn} which finishes the construction of $n_k$, $\varepsilon_k$, $u_k$ and $a$.

Remark that by \eqref{fvdjvjfvfnh} and \eqref{scjkvhsdjovh} we obtain (with a convention $\sum_1^0=0$)
\begin{align}
&\sum_{k=1}^\infty\frac{\varepsilon_k^q}{k}\sum_{j=1}^{n_{k-1}}j^{q/p-1}\overset{L \ref{vjifdjkvfjo}}{\lesssim}
\sum_{k=2}^\infty\frac{\varepsilon_k^q}{k} n_{k-1}^{q/p}\overset{\eqref{scjkvhsdjovh}}{\le}\sum_{k=2}^\infty\frac{c_{k-1}}{k} n_{k-1}^{q/p}\label{wdivkwuv}\\
&\overset{\eqref{fvdjvjfvfnh}}{\le}\sum_{k=2}^\infty\frac{1}{k}\frac{1}{(k-1)n_{k-1}^{q/p}} n_{k-1}^{q/p}=\sum_{k=2}^\infty\frac{1}{k(k-1)}:=B<\infty.\nonumber
\end{align}
Remark that due to \eqref{skjcdwlkjcjdwc}, \eqref{flkkbjbb} and the embedding $\ell^{p,q}\hookrightarrow \ell^{p,r}$ we have
\begin{align}
&\|a\|_{p,r}\le D_{q,r}\|a\|_{p,q}\le D_{q,r}\|\widetilde{a}\|_{p,q}\le D_{q,r}.\label{wdkjcwklcelkck}
\end{align}

Set
\begin{align*}
&z_N=\sum_{k=1}^N k^{-1/q}\ u_k.
\end{align*}
Then $z_N\in X$. Estimate
\begin{align}
&\|z_N\|_{p,q}=\Big\|\sum_{k=1}^N k^{-1/q}\ u_k\Big\|_{p,q}=\Big\|\sum_{k=1}^N k^{-1/q}\ v_k+\sum_{k=1}^N k^{-1/q}\ w_k\Big\|_{p,q}\label{vsdjnhhfioviop}\\
&\overset{\eqref{sdcdiovnhdo}}{\ge} \frac{1}{T}\Big\|\sum_{k=1}^N k^{-1/q}\ v_k\Big\|_{p,q}-\Big\|\sum_{k=1}^N k^{-1/q}\ w_k\Big\|_{p,q}.\nonumber
\end{align}
Clearly we have
\begin{align}
&\Big\|\sum_{k=1}^N k^{-1/q}\ w_k\Big\|_{p,q}\overset{\eqref{dovdovio}}{\le} \sum_{k=1}^N k^{-1/q} T^k\|w_k\|_{p,q}\overset{\text{\eqref{elbkjbpbjp}}}{\le}\sum_{k=1}^N k^{-1/q} T^k\frac{\delta}{(2T)^{k}}\label{pjdfbioedvfd}\\
&\le\sum_{k=1}^N \frac{\delta}{2^{k}}=\delta.\nonumber
\end{align}
Denote $A_k=\{i\in I_k; v_k(i)=0\}$ and define
\begin{align*}
&a_k(j)=(a\chi_{I_k})(j),\ j\in\mathbb{N},\\
&\tilde{v}_k(j)=|v_k(j)|+a_k(j)\chi_{A_k}(j),\ j\in\mathbb{N}.
\end{align*}
where $a$ is the fixed constructed sequence.

Fix now $i\in I_k$, $j\in I_{k+1}$. If $v_k(i)\neq 0$ then
\begin{align*}
&\widetilde{v}_k(i)=|v_k(i)|\ge b_k\overset{\eqref{scjkvhsdjovh}}{\ge}\varepsilon_{k+1}\overset{\eqref{ebjklnbkbnkpn}}{\ge}|v_{k+1}(j)|
\end{align*}
and also
\begin{align*}
&\widetilde{v}_k(i)\ge b_k\overset{\eqref{dsvjdvjvjwvkjv}}{\ge}a(n_k+1)\ge a(j).
\end{align*}
So
\begin{align*}
&\widetilde{v}_k(i)\ge |v_{k+1}(j)|+a(j)\chi_{A_{k+1}}(j)=\widetilde{v}_{k+1}(j).
\end{align*}
If $v_k(i)= 0$ then
\begin{align*}
&\widetilde{v}_k(i)=a(i)\ge a(n_k)\overset{\eqref{scjkvhsdjovh}}{\ge}\varepsilon_{k+1}\overset{\eqref{ebjklnbkbnkpn}}{\ge}|v_{k+1}(j)|
\end{align*}
and also
\begin{align*}
&\widetilde{v}_k(i)\ge a(n_k)\ge a(n_k+1)\ge a(j)
\end{align*}
which gives again
\begin{align*}
&\widetilde{v}_k(i)\ge |v_{k+1}(j)|+a(j)\chi_{A_{k+1}}(j)=\widetilde{v}_{k+1}(j).
\end{align*}

It implies $\tilde{v}_k(i)\ge\tilde{v}_{k+1}(j)$ for $i\in I_k, j\in I_{k+1}$ which yields immediately
\begin{align}
&k^{-1/q}\tilde{v}_k(i)>(k+1)^{-1/q}\tilde{v}_{k+1}(j),\ \ \  i\in I_k, j\in I_{k+1}\label{wdlkcjwlkcj}
\end{align}
and so,
\begin{align}
&\Big(\sum_{k=1}^N k^{-1/q}\ \tilde{v}_k\Big)^*=\Big(\sum_{k=1}^N k^{-1/q}\ \tilde{v}_k^{\diamond}\Big)=\sum_{k=1}^N k^{-1/q}\sum_{j=n_{k-1}+1}^{n_k} \tilde{v}_k^{\diamond}\chi_{\{j\}}.\label{wdlvjvjw}
\end{align}
Since $\supp v_k$ are pairwise disjoint we have
\begin{align*}
&\Big\|\sum_{k=1}^N k^{-1/q}\ v_k\Big\|_{p,q}=\Big\|\sum_{k=1}^N k^{-1/q}\ |v_k|\Big\|_{p,q}\\
&=\Big\|\sum_{k=1}^N k^{-1/q}\ \tilde{v}_k-\sum_{k=1}^N k^{-1/q}\ a_k\chi_{A_k}\Big\|_{p,q}\\
&\overset{\eqref{sdcdiovnhdo}}{\ge} \frac{1}{T}\Big\|\sum_{k=1}^N k^{-1/q}\ \tilde{v}_k\Big\|_{p,q}-\Big\|\sum_{k=1}^N k^{-1/q}\ a_k\Big\|_{p,q}.
\end{align*}
Since $k^{-1/q}\le 1$ and $\supp a_k$ are pairwise disjoint we have
\begin{align*}
&\Big\|\sum_{k=1}^N k^{-1/q}\ a_k\Big\|_{p,q}\le\|a\|_{p,q}
\end{align*}
which concludes

\begin{align}
&\Big\|\sum_{k=1}^N k^{-1/q}\ v_k\Big\|_{p,q}\ge\frac{1}{T}\Big\|\sum_{k=1}^N k^{-1/q}\ \tilde{v}_k\Big\|_{p,q}-\|a\|_{p,q}.\label{foijwriopbjpb}
\end{align}
Further
\begin{align*}
&\Big\|\sum_{k=1}^N k^{-1/q}\ \widetilde{v}_k\Big\|_{p,q}^q=\Big\|\sum_{k=1}^N k^{-1/q}\sum_{j=n_{k-1}+1}^{n_k} \widetilde{v}_k(j)\Big\|_{p,q}^q\\
&\overset{\eqref{wdlkcjwlkcj}}{=}\Big\|\sum_{k=1}^N\sum_{j=n_{k-1}+1}^{n_k} (k^{-1/q}\widetilde{v}_k(j))^\diamond\Big\|_{p,q}^q=
\Big\|\sum_{k=1}^N\sum_{j=n_{k-1}+1}^{n_k}k^{-1/q} (\widetilde{v}_k)^\diamond(j)\Big\|_{p,q}^q\\
&\overset{\eqref{wdlvjvjw}}{=}\sum_{k=1}^N\sum_{j=n_{k-1}+1}^{n_k}j^{q/p-1}k^{-1} ((\widetilde{v}_k)^\diamond(j))^q=\sum_{k=1}^N k^{-1}\sum_{j=n_{k-1}+1}^{n_k}j^{q/p-1} ((\widetilde{v}_k)^\diamond(j))^q\\
&=\sum_{k=1}^N k^{-1}\sum_{j=1}^{n_k -n_{k-1}}(j+n_{k-1})^{q/p-1} ((\widetilde{v}_k)^\diamond(j+n_{k-1}))^q\\
&\overset{\eqref{vjrvjvjvvvj}}{=}\sum_{k=1}^N k^{-1}\sum_{j=1}^{n_k -n_{k-1}}(j+n_{k-1})^{q/p-1} ((\widetilde{v}_k)^*(j))^q\\
&\ge\sum_{k=1}^N k^{-1}\sum_{j=n_{k-1}+1}^{n_k -n_{k-1}}(j+n_{k-1})^{q/p-1} ((\widetilde{v}_k)^*(j))^q.
\end{align*}

Since $\tilde{v}_k^{*}(j)\ge {v}_k^{*}(j)$ we obtain
\begin{align}
&\Big\|\sum_{k=1}^N k^{-1/q}\ \tilde{v}_k\Big\|_{p,q}^q\ge\sum_{k=1}^N k^{-1}\sum_{j=n_{k-1}+1}^{n_k -n_{k-1}} (j+n_{k-1})^{q/p-1} ({v}_k^{*}(j))^{q}\label{sdjvkhsduojbvhiobv}\\
&=\sum_{k=1}^N k^{-1}\sum_{j=n_{k-1}+1}^{n_k -n_{k-1}} \Big(\frac{j+n_{k-1}}{j}\Big)^{q/p-1} j^{q/p-1} ({v}_k^{*}(j))^{q}.\nonumber
\end{align}
Since
\begin{align*}
&1\le\frac{j+n_{k-1}}{j}\le 2\ \ \text{for}\ \ n_{k-1}+1\le j
\end{align*}
we have
\begin{align}
&\min\{1,2^{q/p-1}\}\le \Big(\frac{j+n_{k-1}}{j}\Big)^{q/p-1}\le \max\{1,2^{q/p-1}\}\label{sfiviovhio}
\end{align}
which yields with \eqref{sdjvkhsduojbvhiobv}
\begin{align*}
&\Big\|\sum_{k=1}^N k^{-1/q}\ \tilde{v}_k\Big\|_{p,q}^q\gtrsim \sum_{k=1}^N k^{-1}\sum_{j=n_{k-1}+1}^{n_k -n_{k-1}}  j^{q/p-1} ({v}_k^{*}(j))^{q}\\
&\ge \sum_{k=1}^N k^{-1}\sum_{j=1}^{n_k -n_{k-1}}  j^{q/p-1} ({v}_k^{*}(j))^{q}\nonumber -\sum_{k=1}^N k^{-1}\sum_{j=1}^{n_{k-1}}  j^{q/p-1} ({v}_k^{*}(j))^{q}.
\end{align*}
Clearly,
\begin{align*}
&\sum_{j=1}^{n_k -n_{k-1}}  j^{q/p-1} ({v}_k^{*}(j))^{q}=\|v_k\|_{p,q}^q
\end{align*}
and so
\begin{align}
&\Big\|\sum_{k=1}^N k^{-1/q}\ \tilde{v}_k\Big\|_{p,q}^q\overset{\text{\eqref{ebjklnbkbnkpn}}}{\gtrsim}\sum_{k=1}^N k^{-1}\|v_k\|_{p,q}^{q} -\sum_{k=1}^N k^{-1}\sum_{j=1}^{n_{k-1}}  j^{q/p-1} \varepsilon_k^{q}\label{dvhiovhiohib}\\
&\overset{\text{\eqref{efbklbpmjkpbn}}}{\ge}   \sum_{k=1}^N k^{-1}\Big(\frac{1}{T}-\delta\Big)^q-\sum_{k=1}^N k^{-1}\sum_{j=1}^{n_{k-1}}  j^{q/p-1} \varepsilon_k^{q}\nonumber\\
&\overset{\text{\eqref{wdivkwuv}}}{\gtrsim } \Big(\Big(\frac{1}{T}-\delta\Big)^q\sum_{k=1}^N k^{-1}-B \Big)\gtrsim  A\ln N-B.\nonumber
\end{align}
Now,
\begin{align*}
&\|z_N\|_{p,q}\overset{\text{\eqref{vsdjnhhfioviop}}}{\ge}\frac{1}{T}\Big\|\sum_{k=1}^N k^{-1/q}\ v_k\Big\|_{p,q}-\Big\|\sum_{k=1}^N k^{-1/q}\ w_k\Big\|_{p,q}\\
&\overset{\text{\eqref{pjdfbioedvfd},\eqref{foijwriopbjpb}}}{\ge}\frac{1}{T}\Big( \frac{1}{T}\Big\|\sum_{k=1}^N k^{-1/q}\ \tilde{v}_k\Big\|_{p,q}-\|a\|_{p,q}\Big)-\delta\\
&\overset{\text{\eqref{dvhiovhiohib}}}{\gtrsim}\frac{1}{T^2}(A\ln N-B)^{1/q}-\frac{\|a\|_{p,q}}{T}-\delta.
\end{align*}
It implies
\begin{align}
&\|z_N\|_{p,q}\rightarrow\infty\ \text{for}\ N\rightarrow\infty.\label{efvjjvjpvjj}
\end{align}

It remains to estimate $\|z_N\|_{p,r}$. Clearly
\begin{align}
&\|z_N\|_{p,r}=\Big\|\sum_{k=1}^N k^{-1/q}\ u_k\Big\|_{p,r}=\Big\|\sum_{k=1}^N k^{-1/q}\ v_k+\sum_{k=1}^N k^{-1/q}\ w_k\Big\|_{p,r}\label{jprjvrjrjj}\\
&\le T\Big(\Big\|\sum_{k=1}^N k^{-1/q}\ v_k\Big\|_{p,r}+\Big\|\sum_{k=1}^N k^{-1/q}\ w_k\Big\|_{p,r}\Big).\nonumber\\
&\overset{\text{\eqref{pjdfbioedvfd}}}{\le} T\Big(\Big\|\sum_{k=1}^N k^{-1/q}\ v_k\Big\|_{p,r}+\delta\Big).\nonumber
\end{align}
Further
\begin{align}
&\Big\|\sum_{k=1}^N k^{-1/q}\ v_k\Big\|_{p,r}^r\le \Big\|\sum_{k=1}^N k^{-1/q}\ \tilde{v}_k\Big\|_{p,r}^r=\Big\|\Big(\sum_{k=1}^N k^{-1/q}\ \tilde{v}_k\Big)^*\Big\|_{p,r}^r\label{evjvjjvjvjv}\\
&\overset{\eqref{wdlvjvjw}}{=}\sum_{k=1}^N k^{-r/q}\sum_{j=n_{k-1}+1}^{n_k} j^{r/p-1} (\tilde{v}_k^{\diamond}(j))^{r}\nonumber\\
&=\sum_{k=1}^N k^{-r/q}\sum_{j=n_{k-1}+1}^{2n_{k-1}} j^{r/p-1} (\tilde{v}_k^{\diamond}(j))^{r}\nonumber   +\sum_{k=1}^N k^{-r/q}\sum_{j=2n_{k-1}+1}^{n_k} j^{r/p-1} (\tilde{v}_k^{\diamond}(j))^{r}.\nonumber
\end{align}
First estimate
\begin{align}
&\sum_{k=1}^N k^{-r/q}\sum_{j=n_{k-1}+1}^{2n_{k-1}} j^{r/p-1} (\tilde{v}_k^{\diamond}(j))^{r}\label{doivcoiwviwv}\\
&=\sum_{k=1}^N k^{-r/q}\sum_{j=n_{k-1}+1}^{2n_{k-1}} j^{r/p-1} ({v}_k+a_k\chi_{A_k})^{\diamond}(j))^{r}\nonumber\\
&\overset{\text{\eqref{ebjklnbkbnkpn}}}{\le}\sum_{k=1}^N k^{-r/q}\sum_{j=n_{k-1}+1}^{2n_{k-1}} j^{r/p-1} (\varepsilon_k+a_k\chi_{A_k})^{\diamond}(j))^{r}.\nonumber
\end{align}
Clearly
\begin{align}
&\sum_{j=n_{k-1}+1}^{2n_{k-1}} j^{r/p-1}\lesssim \int_{n_{k-1}}^{2n_{k-1}}x^{r/p-1}dx\lesssim n_{k-1}^{r/p}.\label{edcjvjpvjvjpjv}
\end{align}

Since $a$ is non-increasing sequence and $\varepsilon_k$ is constant on $I_k$ we have $(\varepsilon_k+a_k)^{\diamond}(j)=\varepsilon_k+a_k(j)$ which implies
\begin{align}
&\sum_{k=1}^N k^{-r/q}\sum_{j=n_{k-1}+1}^{2n_{k-1}} j^{r/p-1} (\varepsilon_k+a_k\chi_{A_k})^{\diamond}(j))^{r}\label{ldjvvjejjj}\\
&\le\sum_{k=1}^N k^{-r/q}\sum_{j=n_{k-1}+1}^{2n_{k-1}} j^{r/p-1} (\varepsilon_k+a_k(j))^{r}\nonumber\\
&\lesssim\Big(\sum_{k=1}^N k^{-r/q}\sum_{j=n_{k-1}+1}^{2n_{k-1}} j^{r/p-1} \varepsilon_k^{r}\nonumber\\
&\qquad\qquad+\sum_{k=1}^N k^{-r/q}\sum_{j=n_{k-1}+1}^{2n_{k-1}} j^{r/p-1} a_k^{r}(j)\Big)\nonumber\\
&\overset{\text{\eqref{edcjvjpvjvjpjv}}}{\lesssim}  \Big(\sum_{k=1}^N k^{-r/q}n_{k-1}^{r/p} \varepsilon_k^{r}+\sum_{k=1}^N \sum_{j=n_{k-1}+1}^{2n_{k-1}} j^{r/p-1} a_k^{r}(j)\Big) \nonumber \\
&\overset{\text{\eqref{efbklbnpbkk}}}{\lesssim} \Big(\sum_{k=1}^\infty k^{-r/q}+\|a\|_{p,r}^r\Big)\overset{\eqref{wdkjcwklcelkck}}{<}\infty.\nonumber
\end{align}

Estimate
\begin{align*}
&\sum_{j=2n_{k-1}+1}^{n_k} j^{r/p-1} (\tilde{v}_k^{\diamond}(j))^{r}\overset{\text{\eqref{vjrvjvjvvvj}}}{=}\sum_{j=n_{k-1}+1}^{n_k-n_{k-1}} (j+n_{k-1})^{r/p-1} (\tilde{v}_k^{*}(j))^{r}\\
&\overset{\text{\eqref{sfiviovhio}}}{\lesssim}\sum_{j=n_{k-1}+1}^{n_k-n_{k-1}} j^{r/p-1} (\tilde{v}_k^{*}(j))^{r}\le\sum_{j=1}^{n_k-n_{k-1}} j^{r/p-1} (\tilde{v}_k^{*}(j))^{r}=\|\tilde{v}_k\|_{p,r}^r\\
&\overset{\text{\eqref{dcehoiefroir}}}{\lesssim}D_{q,r}^r\|\tilde{v}_k\|_{p,q}^r{\lesssim}\|{v}_k+a_k\chi_{A_k}\|_{p,q}^r\lesssim\|{v}_k\|_{p,q}^r+\|a_k\|_{p,q}^r\\
&\overset{\text{\eqref{efbklbpmjkpbn}}}{\lesssim} 1+\|a\|_{p,q}^r:=C<\infty.
\end{align*}
Thus
\begin{align}
&\sum_{k=1}^N k^{-r/q}\sum_{j=2n_{k-1}+1}^{n_k} j^{r/p-1} (\tilde{v}_k^{\diamond}(j))^{r}\le C\sum_{k=1}^N k^{-r/q}<\infty.\label{woiceiiu}
\end{align}
Now, \eqref{jprjvrjrjj}, \eqref{evjvjjvjvjv}, \eqref{doivcoiwviwv}, \eqref{ldjvvjejjj} and \eqref{woiceiiu} show that
\begin{align*}
&\|z_N\|_{p,r}\le K<\infty\ \text{far all}\ N
\end{align*}
and this with \eqref{efvjjvjpvjj} finishes the proof.
\end{proof}

Now  we will investigate the case $q=\infty$.
\begin{lem}\label{coeifhiehivh}
Let $X,Y,Z$ be quasi-Banach spaces and $T:X\rightarrow Y$, $S:Y\rightarrow Z$ be linear bounded mappings. Assume that either $T$ is strictly singular. Then the composition $T\circ S:X\rightarrow Z$ is strictly singular.
\end{lem}
\begin{proof}
Let $T\circ S$ is not strictly singular. Then there are infinite dimensional subspace $P\subset X$ and positive constant $c_1,c_2$ such that for all $u\in P$ we have
\begin{align*}
&c_1\|S(T(u))\|_Z\le \|u\|_X\le c_2\|S(T(u)).
\end{align*}
Due to the boundedness of $T$ we obtain for each $u\in P$
\begin{align*}
&\frac{1}{\|T\|}\|T(u)\|_Y\le \|u\|_X\le c_2\|S(T(u))\le c_2\|S\|\ \|T(u)\|_Y
\end{align*}
and so $P, T(P)$ are isomorphic. It contradicts the assumption $T$ is strictly singular.
\end{proof}

\begin{thm}
Let $0< p<\infty$, $0<q<r\le \infty$. Then the embedding $\ell_{p,q}\hookrightarrow\ell_{p,r}$ is strictly singular.
\end{thm}
\begin{proof}
The case $q<\infty$ is proved in Theorem \ref{wdjkchdocvho}. Assume now $q=\infty$. Choose $s, \ q<s<\infty$. Then
\begin{align*}
&\ell^{p,q}\hookrightarrow \ell^{p,s}\hookrightarrow \ell^{p,\infty}.
\end{align*}
Since $\ell^{p,q}\hookrightarrow \ell^{p,s}$ is strictly singular by Theorem \ref{wdjkchdocvho} we have by Lemma \ref{coeifhiehivh} that the embedding $\ell^{p,q}\hookrightarrow \ell^{p,\infty}$ considering as a composition of two embeddings is strictly singular.
\end{proof}

\section{The embedding is not finitely strictly singular}

\begin{defi}
Given $n\in\mathbb{N}$ we define a function on an interval $[0,1]$ by
\begin{align*}
R_n(t)=\sign \sin 2^n \pi t.
\end{align*}
\end{defi}
We can describe $R_n$ by more natural way. Set $I_i=(\frac{i-1}{2^n},\frac{i}{2^n})$, $i=1,2,\dots,2^n$ we have
\begin{align*}
R_n(t)=\sum_{i=1}^{2^n}(-1)^{i+1}\chi_{I_i}(t).
\end{align*}

Let us remind a well-known Khintchine's inequality. A proof can be found in \cite{A}, Theorem 1.4.
\begin{thm}[Khintchine's inequality]\label{eowjeoijoiijfpoi}
Let $0<p<\infty$. Then there are constants $A_p,B_p$ such that for all $N\in\mathbb{N}$ and $a=(a_1,a_2,\dots,a_N)\in \mathbb{R}^N$ we have
\begin{align*}
A_p\Big\|\sum_{i=1}^N a_iR_i\Big\|_{L^p(0,1)}\le \|a\|_{\ell^2}:=\Big(\sum_{i=1}^N a_i^2\Big)^{1/2} \le B_p\Big\|\sum_{i=1}^N a_iR_i\Big\|_{L^p(0,1)}
\end{align*}
\end{thm}

Given $a=(a_1,a_2,\dots,a_{2^n})\in \mathbb{R}^{2^n}$ we can assign to this sequence a function $A(.)$ defined on $[0,1]$ by
\begin{align*}
A(t)=\sum_{i=1}^{2^n} a_i \chi_{I_i}(t).
\end{align*}

\begin{lem}\label{owjvcwvjrwpv}
Let $0<p,q<\infty$. Then for all $a=(a_1,a_2,\dots,a_{2^n})\in \mathbb{R}^{2^n}$ we have
\begin{align*}
\min\Big\{\Big(\frac{p}{q}\Big)^{1/q},2^{\frac{1}{q}-\frac{1}{p}}\Big\}\|a\|_{\ell^{p,q}}\le 2^{\frac{n}{p}}\|A\|_{L^{p,q}(0,1)} \le \|a\|_{\ell^{p,q}}.
\end{align*}
\end{lem}
\begin{proof}Fix $a=(a_1,a_2,\dots,a_{2^n})\in \mathbb{R}^{2^n}$. Consider the function $A$. Then
\begin{align*}
&\|A\|_{L^{p,q}(0,1)}=\Big\|\sum_{i=1}^{2^n} a_i \chi_{I_i}(t)\Big\|_{L^{p,q}(0,1)}=\Big\|\sum_{i=1}^{2^n} a_i^* \chi_{I_i}(t)\Big\|_{L^{p,q}(0,1)}\\
&=\Big(\int_0^1 t^{\frac{q}{p}-1}\Big(\sum_{i=1}^{2^n} a_i^* \chi_{I_i}(t)\Big)^q dt\Big)^{1/q}=
\Big(\sum_{i=1}^{2^n}\int_{I_i} t^{\frac{q}{p}-1}(a_i^*)^q dt\Big)^{1/q}\\
&=\Big(\sum_{i=1}^{2^n}(a_i^*)^q\int_{\frac{i-1}{2^n}}^{\frac{i}{2^n}} t^{\frac{q}{p}-1} dt\Big)^{1/q}.
\end{align*}
Clearly,
\begin{align*}
&\|A\|_{L^{p,q}(0,1)}=\Big(\sum_{i=1}^{2^n}(a_i^*)^q\int_{\frac{i-1}{2^n}}^{\frac{i}{2^n}} t^{\frac{q}{p}-1} dt\Big)^{1/q}
\le\Big(\sum_{i=1}^{2^n}(a_i^*)^q \frac{1}{2^n}\Big(\frac{i}{2^n}\Big)^{\frac{q}{p}-1} \Big)^{1/q}\\
&=\frac{1}{2^{n/p}}\Big(\sum_{i=1}^{2^n}(a_i^*)^q i^{\frac{q}{p}-1}\Big)^{1/q}=\frac{1}{2^{n/p}}\|a\|_{\ell^{p,q}}
\end{align*}
which proves the inequality
\begin{align*}
2^{\frac{n}{p}}\|A\|_{L^{p,q}(0,1)} \le \|a\|_{\ell^{p,q}}.
\end{align*}
Estimate the second inequality. We have
\begin{align*}
&\|A\|_{L^{p,q}(0,1)}=\Big(\sum_{i=1}^{2^n}(a_i^*)^q\int_{\frac{i-1}{2^n}}^{\frac{i}{2^n}} t^{\frac{q}{p}-1} dt\Big)^{1/q}\\
&=\Big((a_1^*)^q\int_{0}^{\frac{1}{2^n}} t^{\frac{q}{p}-1} dt+\sum_{i=2}^{2^n}(a_i^*)^q\int_{\frac{i-1}{2^n}}^{\frac{i}{2^n}} t^{\frac{q}{p}-1} dt\Big)^{1/q}\\
&\ge \Big(\frac{p}{q}\Big(\frac{1}{2^n}\Big)^{\frac{q}{p}}(a_1^*)^q
+\sum_{i=2}^{2^n}(a_i^*)^q \frac{1}{2^n}\Big(\frac{i-1}{2^n}\Big)^{\frac{q}{p}-1} \Big)^{1/q}\\
&=\frac{1}{2^{n/p}}\Big(\frac{p}{q}(a_1^*)^q+\sum_{i=2}^{2^n}(a_i^*)^q \Big(\frac{i-1}{i}\Big)^{\frac{q}{p}-1} i^{\frac{q}{p}-1}\Big)^{1/q}\\
&\ge\frac{1}{2^{n/p}}\Big(\frac{p}{q}(a_1^*)^q+\sum_{i=2}^{2^n}(a_i^*)^q \Big(\frac{1}{2}\Big)^{\frac{q}{p}-1} i^{\frac{q}{p}-1}\Big)^{1/q}\\
&\ge\frac{1}{2^{n/p}}\min\Big\{\Big(\frac{p}{q}\Big)^{1/q},\Big(\frac{1}{2}\Big)^{\frac{1}{p}-\frac{1}{q}}\Big\}\Big((a_1^*)^q+\sum_{i=2}^{2^n}(a_i^*)^q i^{\frac{q}{p}-1}\Big)^{1/q}\\
&=\frac{1}{2^{n/p}}\min\Big\{\Big(\frac{p}{q}\Big)^{1/q},2^{\frac{1}{q}-\frac{1}{p}}\Big\}\Big(\sum_{i=1}^{2^n}(a_i^*)^q i^{\frac{q}{p}-1}\Big)^{1/q}\\
&=\frac{1}{2^{n/p}}\min\Big\{\Big(\frac{p}{q}\Big)^{1/q},2^{\frac{1}{q}-\frac{1}{p}}\Big\}\|a\|_{\ell^{p,q}}
\end{align*}
which proves the inequality
\begin{align*}
\min\Big\{\Big(\frac{p}{q}\Big)^{1/q},2^{\frac{1}{q}-\frac{1}{p}}\Big\}\|a\|_{\ell^{p,q}}\le 2^{\frac{n}{p}}\|A\|_{L^{p,q}(0,1)}.
\end{align*}
\end{proof}
\begin{lem}\label{kovevjrvjr}
Let $0<p<\infty$. Then for all $a=(a_1,a_2,\dots,a_{2^n})\in \mathbb{R}^{2^n}$ we have
\begin{align*}
\|a\|_{\ell^{p,\infty}}=2^{n/p}\|A\|_{L^{p,\infty}(0,1)}.
\end{align*}
\end{lem}
\begin{proof}Fix $a=(a_1,a_2,\dots,a_{2^n})\in \mathbb{R}^{2^n}$. Consider the function $A$. Then
\begin{align*}
&\|A\|_{L^{p,\infty}(0,1)}=\Big\|\sum_{i=1}^{2^n} a_i \chi_{I_i}(t)\Big\|_{L^{p,\infty}(0,1)}=\Big\|\sum_{i=1}^{2^n} a_i^* \chi_{I_i}(t)\Big\|_{L^{p,\infty}(0,1)}\\
&=\sup_{t\in(0,1)}t^{1/p}\sum_{i=1}^{2^n} a_i^* \chi_{I_i}(t)=\max_{i=1,2,\dots,2^n}\Big(\frac{i}{2^n}\Big)^{1/p}a_i^*=2^{-n/p}\max_{i=1,2,\dots,2^n}i^{1/p}a_i^*\\
&=2^{-n/p}\|a\|_{\ell^{p,\infty}}.
\end{align*}
\end{proof}

The next lemma is an easy modification of Khintchine's inequality for Lorentz spaces.
\begin{lem}\label{evjpejvrpopoj}
 Let $0<p<\infty, 0<q\le\infty$. Then there is a constant $C_{p,q}$ such that
for all $N\in\mathbb{N}$ and $a=(a_1,a_2,\dots,a_N)\in \mathbb{R}^N$ we have
\begin{align*}
C_{p,q}^{-1}\Big\|\sum_{i=1}^N a_iR_i\Big\|_{L^{p,q}(0,1)}\le \|a\|_{\ell^2}\le C_{p,q}\Big\|\sum_{i=1}^N a_iR_i\Big\|_{L^{p,q}(0,1)}.
\end{align*}
\end{lem}
\begin{proof} Take $p_1<p<p_2$. Then $L^{p_2}(0,1)\hookrightarrow L^{p,q}(0,1)\hookrightarrow L^{p_1}(0,1)$. Then there is $K>0$ such that for all $u$ we have
\begin{align*}
K^{-1}\|u\|_{L^{p_1}(0,1)}\le \|u\|_{L^{p,q}(0,1)}\le K\|u\|_{L^{p_2}(0,1)}.
\end{align*}
Consider a function $\sum_{i=1}^N a_iR_i$. By Theorem \ref{eowjeoijoiijfpoi}  we have
\begin{align*}
&\frac{A_{p_2}}{K}\Big\|\sum_{i=1}^N a_iR_i\Big\|_{L^{p,q}(0,1)}\le A_{p_2}\Big\|\sum_{i=1}^N a_iR_i\Big\|_{L^{p_2}(0,1)}\le \|a\|_{\ell^2}\\
&\le B_{p_1}\Big\|\sum_{i=1}^N a_iR_i\Big\|_{L^{p_1}(0,1)}\le KB_{p_1}\Big\|\sum_{i=1}^N a_iR_i\Big\|_{L^{p,q}(0,1)}.
\end{align*}
\end{proof}

Consider an embedding $\ell_{p,q}\rightarrow \ell_{p,r}$ with $q<r$ . Given $L_n\subset\subset \ell_{p,q}$ with $\dim L_n=n$ set
$$
b_n(L_n)=\inf_{u\in L_n}\frac{\|u\|_{p,r}}{\|u\|_{p,q}}.
$$

Say that $M\subset\mathbb{N}$ is an interval if for all $i<j<k$ we have $ j\in M$  provided $i,k\in M$. Let $M_1,M_2\subset \mathbb{N}$. Say that $ M_1\prec M_2$ if for all $i\in M_1, j\in M_2$ we have $i<j$.

Let $n\in \mathbb{N}$ and let $1\le i\le 2^n$. Split a set on indices $1,2,\dots, 2^n$ into $2^{i}$ pairwise disjoint intervals $M_k$, $k=1,2,\dots, 2^i$ such that
$\# M_k=2^{n-i}$, $M_k\prec M_l$ provided $k<l$. Define now sequences $r_{i,n}$ by
\begin{align*}
&r_{i,n}(j)=\begin{cases}(-1)^{k+1} &\ j\in M_k\\
                           0          &\ j>2^n.
             \end{cases}
\end{align*}
For given $n$ we can see
\begin{align*}
r_{n,n}&=\overbrace{1,-1,\ \ 1,-1,\ \ 1,-1,\ \ 1,-1,\dots,-1,}^{2^{n-1}}\overbrace{\ \ 1,\dots,\ \ 1, -1,\ \ 1,-1,\ \ 1,-1}^{2^{n-1}},0,0,\dots\\
r_{n-1,n}&=1,\ \ 1,-1,-1,\ \ 1,\ \ 1,- 1,-1,\dots,-1,\ \ 1,\dots,- 1, -1,\ \ 1,\ \ 1,- 1,-1,0,0,\dots\\
r_{n-2,n}&=1,\ \ 1,\ \ 1,\ \ 1,-1,- 1,- 1,-1,\dots,-1,\ \ 1,\dots,\ \  1, \ \ 1,- 1,- 1,- 1,-1,0,0,\dots\\
&\vdots\\
r_{1,n}&=1,\ \ 1,\ \ 1,\ \ 1,\ \ 1,\ \ 1,\ \  1,\ \ 1,\dots,\ \ \ 1, - 1,\dots,-  1, -1,- 1,- 1,- 1,-1,0,0,\dots.
\end{align*}
It is a analogous system to Rademacher functions. We can write also
\begin{align*}
r_{i,n}(j)&=\begin{cases}\sign\sin 2^{i-n}\pi j&\ 1\le j\le 2^n\\
                          0&\ j>2^n.
       \end{cases}
\end{align*}

Remark that the appropriate function to $r_{i,n}$ is $R_{i}$.
\begin{thm}
Let $0<p<\infty, 0<q<r\le \infty$. Then there is $\alpha>0$ such that for all $n\in \mathbb{N}$
\begin{align*}
b_n:=b_n(\id:\ell^{p,q}\hookrightarrow\ell^{p,r})\ge \alpha.
\end{align*}
\end{thm}
\begin{proof}
Consider $\mathcal{R}_n=\spa(r_{1,n},r_{2,n},\dots,r_{n,n})$. Then $\dim(\mathcal{R}_n)=n$. Moreover, the appropriate function to a linear combination $\sum_{i=1}^n a_ir_{i,n}$ is a function $A=\sum_{i=1}^n a_iR_i(.)$.
Clearly, by Lemma \ref{kovevjrvjr} and Lemma \ref{owjvcwvjrwpv}
\begin{align*}
&b_n=\sup_{L_n}\inf_{a\in L_n}\frac{\|a\|_{\ell^{p,r}}}{\|a\|_{\ell^{p,q}}}\ge \inf_{a\in\mathcal{R}_n }\frac{\|a\|_{\ell^{p,r}}}{\|a\|_{\ell^{p,q}}}\ge \inf_{a\in\mathcal{R}_n }\min\Big\{\Big(\frac{p}{q}\Big)^{1/q},2^{\frac{1}{q}-\frac{1}{p}}\Big\}\frac{2^{\frac{n}{p}}\|A\|_{L^{p,r}(0,1)}}{2^{\frac{n}{p}}\|A\|_{L^{p,q}(0,1)}}.
\end{align*}
Using Lemma \ref{evjpejvrpopoj} we obtain
\begin{align*}
&b_n\ge C\inf_{a\in\mathcal{R}_n }\frac{\|A\|_{L^{p,r}(0,1)}}{\|A\|_{L^{p,q}(0,1)}}=C\frac{\|\sum_{i=1}^n a_iR_i\|_{L^{p,r}(0,1)}}{\|\sum_{i=1}^n a_iR_i\|_{L^{p,q}(0,1)}}\\
&\ge C\frac{C_{p,r}^{-1}\|a\|_{\ell^2}}{C_{p,q}\|a\|_{\ell^2}}=\frac{C}{C_{p,r}C_{p,q}}:=\alpha.
\end{align*}
\end{proof}

%

\end{document}